\documentclass{article}
\usepackage{amsmath,amsfonts,graphicx}

\evensidemargin
-0.08in \marginparwidth 0.00in
\newtheorem{Lemma}{Lemma}[section]
\newtheorem{Theorem}{Theorem}[section]

\newtheorem{Proposition}{Proposition}[section]

\begin{document}

\title{\bf  Asymptotic stability of a composite wave of
two viscous shock waves for a one-dimensional system of
non-viscous and heat-conductive ideal gas}

\date{} \maketitle
\noindent${\bf {Lili \ Fan^\ast}} \ and \ {\bf {Akitaka \ Matsumura^{\ast\ast}}}$\\
 {\footnotesize{\noindent $\ast$ School of Mathematics and Computer Science,\\
 Wuhan Polytechnic University, Wuhan, China\\
E-mail: fll810@live.cn\\
 \noindent $\ast\ast$ Department of Pure and Applied Mathematics,
 Graduate School of Information \\
 Science and Technology,
 Osaka University, Suita,Osaka, Japen.\\
E-mail: akitaka@ist,osaka-u.ac.jp}}\\

\noindent{\bf{Abstract:} }
This paper is concerned with the asymptotic stability of a composite wave
consisting of two viscous shock waves to the Cauchy problem for a
one-dimensional system of heat-conductive ideal gas without viscosity.
We extend the results by Huang-Matsumura \cite{Huang-Matsumura}
where they treated the
equation of viscous and heat-conductive ideal gas. That is, even for
the non-viscous and heat-conductive case, we show
that if the strengths of the viscous shock waves and
the initial perturbation are suitably
small, the unique global solution in time exists and asymptotically tends
toward the corresponding composite wave whose spacial shifts
of two viscous shock waves are uniquely determined by the initial
perturbation.\\

\section{Introduction}
The one-dimensional motion
of the heat-conductive ideal gas without the viscosity
is described in the Lagrangian mass coordinates
by the system
\begin{eqnarray}\label{1.1}
\left\{
\begin{array}{ll}
v_t-u_x=0,\\[2mm]
u_t+p_x=0,\\[2mm]
(e+\frac{u^{2}}{2})_{t}+(pu)_x=\kappa\big(\dfrac{ \theta_x}{v}\big)_x,
\end{array}
\right.
\end{eqnarray}
where $x\in {\bf
R}$, $t>0$ and $v >0$, $u$, $\theta>0$, $e >0$ and $p$ are
the specific volume, fluid velocity,
internal energy per unit mass,
absolute temperature and pressure respectively, while
the coefficient of heat conduction $\kappa$ is assumed to be a positive constant.
Here we study the ideal and polytropic fluids
  so that $p$ and $e$ are given by the state equations
\begin{equation}\label{1.2}
p=\frac{R\theta}{v},\quad  e=\frac{R}{\gamma-1}\theta
\end{equation}
where $\gamma>1$ is the adiabatic exponent and $R>0$ is the gas constant. In the present paper,
we investigate
the Cauchy problem for (1.1) with the initial data
\begin{equation}\label{1.3}
(v,u,\theta)(x,0)=(v_0,u_0,\theta_0)(x),
\quad x\in {\bf R}, \quad \
\lim\limits_{x\to \pm\infty}(v_0,u_0,\theta_0)(x)=(v_\pm,u_\pm,\theta_\pm),
\end{equation}
where the far field states $v_\pm>0$, $\theta_\pm>0$ and $u_\pm\in {\bf
R}$ are given constants. Then, we are interested in the global solutions
in time of the Cauchy problem (1.1),(1.3) and their large-time behaviors
in the relations with the far field states
$(v_\pm,u_\pm,\theta_\pm)$. In this aspect, there have been many works
in the case where the ideal gas is both viscous and heat-conductive
(cf. [1-5], [8-10], [12,13], [15], [17]):
\begin{eqnarray}\label{1.4}
\left\{
\begin{array}{ll}
v_t-u_x=0,\\[2mm]
u_t+p_x= \mu(\dfrac {u_{x}}v)_{x},\qquad x\in {\bf R}, \ t>0,\\[2mm]
(e+\frac{u^{2}}{2})_{t}+(pu)_x=
\kappa\big(\dfrac{ \theta_x}{v}\big)_x+\mu
\big(\dfrac{uu_{x}}{v}\big)_{x},
\end{array}
\right.
\end{eqnarray}
where the positive constant $\mu$ is the coefficient of viscosity.
In these previous works, it has been known that
the large-time behaviors are well characterized by
the solutions (``Riemann solutions'') of the corresponding Riemann
problem for the hyperbolic part of (\ref{1.1}) (Euler system):
\begin{equation}\label{1.5}
\left\{
\begin{array}{ll}
v_t-u_x=0,\\[2mm]
u_t+p_x=0,\qquad x\in {\bf R}, \ t>0,\\[2mm]
(e+\frac{u^{2}}{2})_{t}+(pu)_x=0,
\end{array}
\right.
\end{equation}
with the Riemann initial data
\begin{equation}\label{1.6}
(v,u,\theta)(x,0)= (v^R_0,u^R_0,\theta^R_0)(x):=
\left\{
\begin{array}{l}
(v_-,u_-,\theta_-), \ x<0,\\[2mm]
(v_+,u_+,\theta_+), \ x>0.
\end{array}
\right.
\end{equation}
Particularly, in the case where the Riemann solution consists of
two shock waves, Huang-Matsumura (\cite{Huang-Matsumura})
investigated the Cauchy problem
(1.4),(1.3) around a linear combination of the corresponding
two viscous shock waves and showed that
if the strengths of the shock waves and initial perturbations are
small, the unique global solution in time exists and asymptotically tends
toward the corresponding composite wave whose spacial shifts
of two viscous shock waves are uniquely determined by the initial
perturbation. The aim of this article is to extend this stability
result to the case where the ideal gas is even non-viscous and heat-conductive.
Our arguments basically follow the arguments in \cite{Huang-Matsumura},
where they
first construct a diffusion wave to reduce the original system (\ref{1.1})
to an integrated system and then apply an elementary
energy method to establish the a priori estimates with the aid of
the monotone property of the viscous shock waves.
However, since our system (\ref{1.1}) as a hyperbolic-parabolic system
is of less dissipative than (\ref{1.4}), we need more
subtle estimates to recover the
regularity and dissipativity for the
components of the hyperbolic part, that is, the
density and fluid velocity. Since our system (\ref{1.1}) still
satisfies a condition of dissipativity (``Kawashima-Shizuta
condition'' (cf. \cite{kawashima}, \cite{kawashima-thesis}),
we shall overcome this difficulty
by manipulating new several energy estimates which correspond
to a dissipative mechanism ensured by the Kawashima-Shizuta
theory, and also looking for the
perturbed solution for the integrated system of (\ref{1.1}) in
$C(H^3)$ to control the nonlinearity of the hyperbolic part, instead of usual $C(H^2)$ in the previous papers.

Finally, we make a remark that in the case where
the Riemann solution consists of
rarefaction waves, Murakami \cite{Murakami} investigated the
asymptotics toward the rarefaction waves
for the non-viscous and heat-conductive
system (1.1), and the case where
the Riemann solution consists of
contact discontinuity and rarefaction waves
will be discussed in a forthcoming paper by
Lili and Murakami.

The rests of the paper are organized as follows:
in Section 2, we  introduce the properties of shock profiles,
 construct the diffusion wave  and state the main results.
In Section 3, we show the way to prove the main theorem
by combining the local existence and a priori estimates.
Finally, in Section 4, the a priori estimates are
established by elementary energy method.

\medskip

\noindent {\bf Notations:}
Throughout this paper, without
any ambiguity, we denote a generic positive constant by $C$ or $c$. If the
dependence needs to be explicitly pointed out, then the notations
$C(\cdot)$ or $C(\cdot,\cdot)$, etc, are used.
For function spaces, $L^{p}(1\leq p \leq \infty)$ denotes
the usual Lebesgue space on
${\bf R}$ with its norm:
\begin{eqnarray*}
&&\|f\|_{L^{p}} = \left(\int_{{\bf R}}
{|f(x)|^{p}}\,dx\right)^{\frac{1}{p}},\quad (1\leq p <\infty),\\[1.5mm]
&&\|f\|_{L^{\infty}} ={\rm ess.}\sup\limits_{\!\!\!\!\!\!\!\!\!\!x\in {\bf R}}
|f(x)|,\quad (p =\infty).
\end{eqnarray*}
$H^k$ denotes the usual $k$-th order Sobolev space with its norm
\begin{eqnarray*}
\|f\|_k: = \bigg(\sum_{j=0}^{k}\|\partial_x^j f\|^2_{L^{2}}\bigg)^{\frac{1}{2}},
\end{eqnarray*}
and when $k=0$, we note $H^0=L^2$ and write its norm $\|\cdot\|$ for simplicity.

\section{Viscous shock waves, diffusion wave and main theorems}
\setcounter{equation}{0}

In this section, we firstly construct the two desired viscous shock waves
for (1.1)
which correspond to the two shock waves making up the Riemann solution of
(1.5). To do that, we recall the Riemann problem (1.5),(1.6) under
consideration.
Setting $z={}^t(v,u,\theta)$, we can rewrite the Euler system (\ref{1.1})
for smooth solution $z$ as in the form
$$
z_t+A(z)z_x=0
$$
where
\begin{eqnarray*}
A(z)=
\left(
\begin{array}{ccc}
0&-1&0\\[2mm]
-\frac{p}{v}&0&\frac{R}{v}\\[2mm]
0&\frac{\gamma-1}{R}p&0
\end{array}
\right).
\end{eqnarray*}
It is known that the Jacobi matrix $A(z)$ has  the three eigenvalues:
$\lambda_1=-\sqrt{\gamma p/v}<0,\lambda_2=0,\lambda_3=-\lambda_1>0$
where the second characteristic field is linear degenerate and the
others are genuinely nonlinear.
In this paper, we consider the situation where
the Riemann solution consists of two shock waves,
that is, there exists an intermediate
state $(v_m,u_m,\theta_m)$ such that $(v_{-},u_{-},\theta_{-})$
connects with $(v_m,u_m,\theta_m)$ by the 1-shock wave with the shock speed
$s_1<0$ and $(v_m,u_m,\theta_m)$ connects with $(v_+,u_+,\theta_+)$
by the 3-shock wave with the shock speed $s_3>0$. Here
the shock speeds $s_1$ and $s_3$ are determined by Rankine-Hugoniot condition
and satisfy entropy conditions
$$
\lambda_1(v_-,\theta_-)>s_1>\lambda_1(v_m,\theta_m),~
\lambda_3(v_m,\theta_m)>s_3>\lambda_3(v_+,\theta_+).
$$
By the standard arguments (e.g.\cite{Smoller}), for each
$(v_{-},u_{-},\theta_{-})$ we can see our
situation takes place provided $(v_+,u_+,\theta_+)$
is located on a quarter of a curved surface in a small neighborhood of $(v_{-},u_{-},\theta_{-})$.
In what follows, the neighborhood of $(v_{-},u_{-},\theta_{-})$ denoted by
$\Omega_-$ is given by
$$
\Omega_- = \{(v,u,\theta) |\ |(v-v_-,u-u_-,\theta-\theta_-)| \le \bar{\delta} \}
$$
where $\bar{\delta}$ is a positive constant depending only on $(v_{-},u_{-},\theta_{-})$.
To describe the
strengths of the shock waves for later use, we set
\begin{eqnarray}\label{2.1}
 &&\delta_1=|v_m-v_-|+|u_m-u_-|+|\theta_m-\theta_-|,\\[2mm]
&&\delta_3=|v_m-v_+|+|u_m-u_+|+|\theta_m-\theta_+|,\nonumber
\end{eqnarray}
and $\delta=\min\{\delta_1,\delta_3\}$.
In our situation, for $(v_+,u_+,\theta_+) \in \Omega_-$,
we note that it holds
\begin{equation}\label{2.2}
\delta_1+\delta_3\leq C|(v_+-v_-,u_+-u_-,\theta_+-\theta_-)|,
\end{equation}
where $C$ is a positive constant depending only on $(v_-,u_-,\theta_-)$.
Then, if it also holds
\begin{equation}\label{2.3}
\delta_1+\delta_3\leq C\delta, \quad  \delta_1+\delta_3\rightarrow0
\end{equation}
for a positive constant $C$, we call the strengths of the shock waves
$\delta_1,\delta_3$  ``small with same order".
We always assume (\ref{2.3}) in what follows.

Now we construct the viscous 1-shock wave of (\ref{1.1})
with the form $(V_1,U_1,\Theta_1)(x-s_1t)$ which corresponds
to the 1-shock wave with the shock speed $s_1<0$ connecting the far field states
$(v_-,u_-,\theta_-)$ and $(v_m,u_m,\theta_m)$.
Substituting $(v,u,\theta)= (V_1,U_1,\Theta_1)(\xi), \xi=x-s_1t$ into (1.1),
it is determined by
\begin{eqnarray}\label{2.4}
\left\{
\begin{array}{ll}
-s_1V'_1-U'_1=0,\\[2mm]
-s_1U'_1+P'_1=0,\\[2mm]
-s_1\bigg(\frac{R}{\gamma-1}\Theta_1+\frac{U^2_1}{2}\bigg)'+(P_1U_1)'=\bigg(\dfrac{\kappa \Theta'_1}{V_1}\bigg)',
\qquad \xi \in {\bf R},\\[2mm]
(V_1,U_1,\Theta_1)(-\infty)=(v_-,u_-,\theta_-),\\[2mm]
(V_1,U_1,\Theta_1)(+\infty)=(v_m,u_m,\theta_m),
\end{array}
\right.
\end{eqnarray}
where $P_1=p(V_1,\Theta_1),p_{\pm}=p(v_\pm,\theta_\pm), e_\pm=e(\theta_\pm)$,
$p_m=p(v_m,\theta_m)$, and $e_m=e(\theta_m)$.
In order for the solution of (2.4) to exist, we can see that
the shock speed $s_1$ and
the far field states must satisfy the Rankine-Hugoniot condition
\begin{eqnarray}\label{2.5}
\left\{
\begin{array}{ll}
-s_1(v_m-v_-)-(u_m-u_-)=0,\\[2mm]
-s_1(u_m-u_-)+(p_m-p_-)=0,\\[2mm]
-s_1\big((e_m+\frac{1}{2}u_m^2)-(e_-+\frac{1}{2}u_-^2)\big)+(p_mu_m-p_-u_-)=0
\end{array}
\right.
\end{eqnarray}
and the entropy condition
\begin{eqnarray}
\lambda_1(v_-,\theta_-)= -\sqrt{\gamma p_-/v_-}>s_1>-\sqrt{\gamma p_m/v_m}=\lambda_1(v_m,\theta_m)
\end{eqnarray}
which is equivalent to $u_->u_m$.
It is noted that the R-H condition (2.5) and entropy condition (2.6) imply
\begin{eqnarray}
s^2_1=\frac{\gamma p_-}{v_m}\bigg(1-\frac{d_-}{1+d_-}\bigg),\quad
\theta_m=\theta_-\bigg(1-\frac{v_-+v_m}{v_-}\frac{d_-}{1+d_-}\bigg),\quad
v_- > v_m
\end{eqnarray}
where $d_-=\frac{\gamma-1}{2}\frac{v_m-v_-}{v_m}<0$.
Here we may assume $\bar{\delta}$, the size of $\Omega_-$, suitably small
to assure that $|d_-|<1$.
Then integrating $(\ref{2.4})$ with respect to $\xi$ under the R-H condition,
we have
\begin{eqnarray}\label{2.6}
\left\{
\begin{array}{ll}
U_1-u_m=-s_1(V_1-v_m),\\[2mm]
P_1-p_m=s_1(U_1-u_m)=-s^2_1(V_1-v_m),\\[2mm]
\dfrac{\kappa \Theta'_1}{V_1}=-s_1\bigg(\frac{R}{\gamma-1}(\Theta_1-\theta_m)+
\frac{1}{2}(U^2_1-u_m^2)\bigg)+(P_1U_1-p_mu_m).
\end{array}
\right.
\end{eqnarray}
Here note that $(v_m, u_m, \theta_m, p_m)$  in (2.8) can be replaced by $(v_-, u_-, \theta_-, p_-)$.
It follows from $(\ref{2.6})_1$ and $(\ref{2.6})_2$ that
\begin{eqnarray}\label{2.100}
\left\{
\begin{array}{ll}
U_1=u_m-s_1(V_1-v_m),\\[2mm]
\Theta_1=\theta_m + \frac{1}{R}\big(-s^2_1(V_1-v_m)^2+(p_m-s^2_1v_m)(V_1-v_m)\big).
\end{array}
\right.
\end{eqnarray}
Hence, once $V_1$ is determined, $U_1$ and $\Theta_1$ are given by
$(\ref{2.100})$. Substituting $(\ref{2.100})$ into $(\ref{2.6})_3$, we have
the equation for $V_1$:
\begin{eqnarray}\label{2.7}
\left\{
\begin{array}{ll}
\kappa V'_1=H_1(V_1),\qquad \xi \in {\bf R},\\[2mm]
V_1(-\infty)=v_-, \quad  V_1(+\infty)=v_m,
\end{array}
\right.
\end{eqnarray}
where
\begin{eqnarray}\label{2.8}
H_1(V_1)=\frac{s_1RV_1\bigg(\frac{\gamma+1}{2}s^2_1(V_1-v_m)^2-(\gamma p_m-s_1^2v_m)(V_1-v_m)\bigg)}
{(\gamma-1)\big(p_m-s^2_1v_m-2s^2_1(V_1-v_m)\big)},
\end{eqnarray}
and note that $v_m$ in $(\ref{2.8})$ can be replaced by $v_-$.
By the entropy condition (2.6) and fact that $s_1^2 \to p_-/v_-\ as\ v_m \to v_-$,
if we take $\bar{\delta}$ properly small if needed, we can easily see that
\begin{equation}
p_m< s^2_1v_m <\gamma p_m, \quad p_-< \gamma p_- < s^2_1v_-.
\end{equation}
Then it follows from (2.11) and (2.12) that
$$
H_1(v_-)=H_1(v_m)=0,\qquad H_1(v)<0\ \ (v_m<v<v_-),
$$
and
\begin{eqnarray*}
H'_1(v_-)=-\frac{s_1Rv_-(\gamma p_--s_1^2v_-)}{(\gamma-1)(p_1-s^2_1v_-)}>0,\quad
H'_1(v_m)=-\frac{s_1Rv_m(\gamma p_m-s_1^2v_m)}{(\gamma-1)(p_m-s^2_1v_m)}<0.
\end{eqnarray*}
Therefore, due to standard theory of ordinary differential equations,
we can construct the solution $V_1$ of (\ref{2.7}) which is unique up to shift
of $\xi$, monotonically decreasing and exponentially tends toward the far field states
$v_-$ and $v_m$.

Similarly,
the viscous 3-shock wave of (\ref{1.1})
with the form $(V_3,U_3,\Theta_3)(\xi), \xi=x-s_3t,$ which corresponds
to the 3-shock wave with the shock speed $s_3>0$ connecting the far field states
$(v_m,u_m,\theta_m)$ and $(v_+,u_+,\theta_+)$ is constructed by
\begin{eqnarray}\label{2.9}
\left\{
\begin{array}{ll}
-s_3V'_3-U'_3=0,\\[2mm]
-s_3U'_3+P'_3=0,\\[2mm]
-s_3\bigg(\frac{R}{\gamma-1}\Theta_3+\frac{U^2_3}{2}\bigg)'
+(P_3U_3)'=(\dfrac{\kappa \Theta'_3}{V_3})',
\qquad \xi \in {\bf R}\\[1.5mm]
(V_3,U_3,\Theta_3)(-\infty)=(v_m,u_m,\theta_m),\\[2mm]
(V_3,U_3,\Theta_3)(-\infty)=(v_+,u_+,\theta_+).
\end{array}
\right.
\end{eqnarray}
under the R-H condition
\begin{eqnarray}\label{2.10}
\left\{
\begin{array}{ll}
-s_3(v_m-v_+)-(u_m-u_+)=0,\\[2mm]
-s_3(u_m-u_+)+(p_m-p_+)=0,\\[2mm]
-s_3\big((e_m+\frac{1}{2}u_m^2)-(e_++\frac{1}{2}u_+^2)\big)+(p_mu_m-p_+u_+)=0,
\end{array}
\right.
\end{eqnarray}
and the entropy condition
\begin{eqnarray*}
\lambda_3(v_m,\theta_m)>s_3>\lambda_3(v_+,\theta_+).
\end{eqnarray*}
Since the remaining arguments on the solution of (2.13) are the same as above
for the 1-viscous shock wave, we omit the details.

Now we are ready to list up the properties of the solutions of (2.4) and (2.13) which
will be used later. Since the computation is standard, we omit the proof.
\begin{Proposition}For any fixed  $(v_-,u_-,\theta_-)$, suppose $(v_+,u_+,\theta_+)\in \Omega_-$
and the Riemann solution of {\rm (1.5)} consist of two shock waves whose strengths
satisfy {\rm (\ref{2.3})}.
Then the problems {\rm (\ref{2.4})} and {\rm (\ref{2.9})} have the smooth solutions
$(V_1,U_1,\Theta_1)(x-s_1t)$ and $(V_3,U_3,\Theta_3)(x-s_3t)$ which are
unique up to spatial shift and satisfy the following:\\
 $(1)$ $U_{ix}(x-s_it)=-s_iV_{ix}(x-s_it)<0$,
 \quad  $|\Theta_{ix}(x-s_it)|\leq C|V_{ix}(x-s_it)|,\quad x\in {\bf R},\ t\geq0$,\quad $i=1,3$;\\[1.5mm]
 $(2)$ There exist some positive constants $c$ and $C$ such that for $i=1,3,$
\begin{eqnarray}\label{2.12}
&&|(V_1-v_m,U_1-u_m,\Theta_1-\theta_m)(x-s_1t)|\leq C\delta_1e^{-c\delta_1|x-s_1t|},
\quad x>s_1t,\ t\geq0,\nonumber\\[2mm]
&&|(V_3-v_m,U_3-u_m,\Theta_3-\theta_m)(x-s_3t)|\leq C\delta_3e^{-c\delta_3|x-s_3t|},
\quad x<s_3t,\ t\geq0,\nonumber\\[2mm]
&& |(V_{ix},U_{ix},\Theta_{ix},V_{ixx},U_{ixx},\Theta_{ixx})(x-s_it)|\leq C\delta^2_ie^{-c\delta_i|x-s_it|},
\quad x\in {\bf R},\ t\geq0,\\[2mm]
&& |s_i^2-\frac{\gamma p_m}{v_m}|\leq C\delta_i.\nonumber
\end{eqnarray}
\end{Proposition}

Now, we turn to the initial value problem (1.1),(1.3) under the situation that
the corresponding Riemann solution consists of two shock waves.
Set $E=\theta+\frac{\gamma-1}{2R}u^2$ and $m(x,t)={}^t(v,u,E)(x,t)$.
Take a pair of viscous shock waves $(V_1,U_1,\Theta_1)(x-s_1t)$ and $(V_3,U_3,\Theta_3)(x-s_3t)$, and
define its composite wave $\overline{m}(x,t)={}^t(\overline{v},\overline{u},\overline{E})(x,t)$ by
\begin{eqnarray}
&&\overline{v}=V_1(x-s_1t)+V_3(x-s_3t)-v_m,\nonumber\\[2mm]
&&\overline{u}=U_1(x-s_1t)+U_3(x-s_3t)-u_m,\nonumber\\[2mm]
&&\overline{E}=E_1(x-s_1t)+E_3(x-s_3t)-E_m,\\[2mm]
&&\overline{\theta}= \overline{E} - \frac{\gamma-1}{2R}{\overline{u}}^2,\nonumber
\end{eqnarray}
where $E_i=\Theta_i+\frac{\gamma-1}{2R}U_i^2$ and $ E_m=\theta_m+\frac{\gamma-1}{2R}u_m^2$.
Here we note that it holds
\begin{equation*}
\overline{\theta} =  \Theta_1(x-s_1t)+\Theta_3(x-s_3t)-\theta_m -
\frac{\gamma-1}{2R}(U_1(x-s_1t)-u_m)(U_3(x-s_3t)-u_m)
\end{equation*}
and Proposition 2.1 and (2.3) give
$$
|(U_1(x-s_1t)-u_m)(U_3(x-s_3t)-u_m)| \le C\delta^2e^{-c\delta(|x|+t)}
$$
which implies that $\overline{\theta}(x,t)$ and $ \Theta_1(x-s_1t)+\Theta_3(x-s_3t)-\theta_m$
are asymptotically equivalent as the time goes to infinity.
Then, we consider the initial value problem (1.1),(1.3) in a small neighborhood of
 $\overline{m}$. As the previous papers show, if the initial mass
$\int (m(x,0)-\overline{m}(x,0))dx$ is zero, we can expect
the solution $m$ tends toward the composite wave $\bar{m}$.
In the case the initial mass
$\int (m(x,0)-\overline{m}(x,0))dx$ is not zero, following the arguments in Liu [10],
the asymptotic state is expected to be a spatially shifted composite wave
$\overline{m}_{\beta_1,\beta_3}={}^t(\overline{v}_{\beta_1,\beta_3},\overline{u}_{\beta_1,\beta_3},
\overline{E}_{\beta_1,\beta_3})$ given by
\begin{eqnarray*}
&&\overline{v}_{\beta_1,\beta_3}=V_1(x-s_1t+\beta_1)+V_3(x-s_3t+\beta_3)-v_m,\nonumber\\[2mm]
&&\overline{u}_{\beta_1,\beta_3}=U_1(x-s_1t+\beta_1)+U_3(x-s_3t+\beta_3)-u_m,\\[2mm]
&&\overline{E}_{\beta_1,\beta_3}=E_1(x-s_1t+\beta_1)+E_3(x-s_3t\beta_3)-E_m, \\[2mm]
&&\overline{\theta}_{\beta_1,\beta_3}= \overline{E}_{\beta_1,\beta_3}
- \frac{\gamma-1}{2R}(\overline{u}_{\beta_1,\beta_3})^2,
\end{eqnarray*}
and the spatial shifts $\beta_1$ and $\beta_3$ are determined
by the following decomposition
\begin{eqnarray}\label{2.13}
\int (m(x,0)-\overline{m}(x,0))\,dx
=\beta_1r_1+\beta_2r_2+\beta_3r_3,
\end{eqnarray}
where $r_1:={}^t(v_m-v_-,u_m-u_-,E_m-E_-)$, $r_3:={}^t(v_+-v_m,u_+-u_m,E_+-E_m)$,
and especially $r_2={}^t(R/p_m,0,1)$ which is the right eigenvector corresponding to
the second eigenvalue $\lambda_2=0$ for $A(z_m)$.

Note that since the vectors $r_1,r_2$, and $r_3$ are linearly independent, the constants
$\beta_1,\beta_2$ and $\beta_3$ are uniquely determined by the initial data.
Since (2.17) implies
\begin{eqnarray}
\int (m(x,0)-{\overline{m}}_{\beta_1,\beta_3}(x,0))\,dx
=\beta_2r_2,
\end{eqnarray}
we can again expect that $m_{\beta_1,\beta_3}$ is the asymptotic state provided
$\beta_2=0$. If $\beta_2$ is not zero, to eliminate the excessive mass $\beta_2r_2$
in the second characteristic field which is linearly degenerate, we introduce a smooth
linear diffusion wave $z^D(x,t)=(v^D,u^D,\theta^D)(x,t)$
($m^D=(v^D,u^D,E^D),\ E^D= \theta^D+\frac{\gamma-1}{2R}(u^D)^2$)
around the constant state $z_m=(v_m,u_m,\theta_m)$ as
in [2]. Here we call $z^D$  a ``diffusion wave'' corresponding to the
second characteristic field when $z^D$
approximately satisfies (1.1) as in the form
\begin{eqnarray}
\left\{
\begin{array}{ll}
v^D_t-u^D_x=0,\\[2mm]
u^D_t+p^D_x= (R^D_1)_x,\\[2mm]
E^D_{t}+(p^Du^D)_x=\kappa\big(\dfrac{ \theta^D_x}{v^D}\big)_x + (R^D_2)_x,
\end{array}
\right.
\end{eqnarray}
where the error terms $R^D_1$ and $R^D_2$ have good decay estimates enough for the {\it a priori} estimates,
and it also satisfies
\begin{equation}
\int (m^D(x,t)-m_m) \,dx= \beta_2 r_2,\quad t \ge 0,
\end{equation}
where $m_m=(v_m,u_m,E_m)$ and the estimates like the linear heat kernel
\begin{equation}
|(z^D(x,t)-z_m)|\le
\frac{C|\beta_2|}{(1+t)^{\frac{1}{2}}}e^{-\frac{cx^2}{1+t}},\quad \ x \in {\bf R},\ t\ge 0,
\end{equation}
for some positive constants $c$ and $C$. We will show how to concretely construct the diffusion
wave $z^D$ later after stating our main theorems.
Once the diffusion wave $z^D$ is defined, we introduce a new
asymptotic state in higher level
$Z={}^t(V,U,\Theta)$
($M={}^t(V,U,\Theta+\frac{\gamma-1}{2R}U^2)$) which is given by $M={\overline{m}}_{\beta_1,\beta_3}+m^D-m_m$,
that is,
\begin{eqnarray}\label{2.18}
V&=&V_1(x-s_1t+\beta_1)+V_3(x-s_3t+\beta_3)-v_m+ (v^D(x,t)-v_m),\nonumber\\[2mm]
U&=&U_1(x-s_1t+\beta_1)+U_3(x-s_3t+\beta_3)-u_m+ (u^D(x,t)-u_m),\\[2mm]
\Theta&+&\frac{\gamma-1}{2R}U^2=E_1(x-s_1t+\beta_1)+E_3(x-s_3t+\beta_3)-E_m+
(E^D(x,t)-E_m).\nonumber
\end{eqnarray}
Then it follows from (2.18) and (2.20) that
\begin{eqnarray}
\int (m(x,0)- M(x,0))\,dx
=0.
\end{eqnarray}
which enable us to apply the anti-derivative method for our problem as in [2] and
to expect the solution $m(x,t)$ tends toward the asymptotic state $M(x,t)$,
eventually ${\overline{m}}_{\beta_1,\beta_3}(x,t)$, as
$t \to \infty$.

We now are ready to state to our main results.
Fixing  $(v_-,u_-,\theta_-)$, we assume that $(v_+,u_+,\theta_+)\in \Omega_-$ and
the Riemann solution of (\ref{1.5}),(1.6) consists of two shock waves.
We consider the asymptotic stability of a composite wave of
two viscous shock waves $(\bar{v}, \bar{u},\bar{\theta})$ defined by (2.16), and
suppose the initial data satisfy
\begin{equation}\label{2.25}
v_0-\bar{v}(0,\cdot),\ u_0-\bar{u}(0,\cdot),\ \theta_0-\bar{\theta}(0,\cdot) \in H^2\cap L^1.
\end{equation}
Then $(V,U,\Theta)$ in (2.22) is well defined where $\beta_1$ and $\beta_3$ are uniquely determined in (\ref{2.13}).
Furthermore, we set
\begin{eqnarray}
\Phi_0(x) &=& \int_{-\infty}^x{(v_0(y)-V(y,0))}\,dy,\nonumber\\[2mm]
\Psi_0(x) &=& \int_{-\infty}^x{(u_0(y)-U(y,0))}\,dy,\\[2mm]
\overline{W}_0(x)&=& \int_{-\infty}^x{\left((\frac{R}{\gamma-1}\theta_0+\frac{u^2_0}{2})(y)-
(\frac{R}{\gamma-1}\Theta+\frac{U^2}{2})(y,0)\right)}\,dy,\nonumber
\end{eqnarray}
and assume
\begin{equation}\label{2.27}
(\Phi_0,\Psi_0,\overline{W}_0)\in L^2,
\end{equation}
which implies $(\Phi_0,\Psi_0,\overline{W}_0)\in H^3$, due to (2.24).
Set
\begin{eqnarray*}
I(v_0,u_0,\theta_0) = \|(v_0-V(\cdot,0), u_0-U(\cdot,0), \theta_0-\Theta(\cdot,0))\|_{H^1\cap L^1}
+\|(\Phi_0,\Psi_0,\overline{W}_0)\|_{L^2}.
\end{eqnarray*}

Then, our main result is the following.
\begin{Theorem}
For any fixed $(v_-,u_-,\theta_-)$, we suppose that $(v_+,u_+,\theta_+)\in \Omega_-$ and
the Riemann solution of {\rm (\ref{1.5}), (\ref{1.6})} consists of two shock waves
whose strengths satisfy {\rm (\ref{2.3})}. Further
assume that the initial data satisfy {\rm (\ref{2.25})} and {\rm (\ref{2.27})}, and
$1\le \gamma< 3$.
Then there
exist positive constants $\delta_0$ and $\epsilon_0$ such
that if $|(v_+-v_-,u_+-u_-,\theta_+-\theta_-)| \le \delta_0$ and
$I(v_0,u_0,\theta_0) \le \epsilon_0$,
then the Cauchy problem {\rm (\ref{1.1})},{\rm (\ref{1.3})} admits a unique
global solution in time $(v,u,\theta)$ satisfying
\begin{eqnarray*}
&&(v-V,u-U,\theta-\Theta)\in C^0([0,\infty);H^2), \\
&&(v-V, u-U) \in L^2(0,\infty;H^2),\quad
\theta-\Theta\in L^2(0,\infty;H^3),
\end{eqnarray*}
and the asymptotic behavior
\begin{eqnarray}\label{2.30}
\lim\limits_{t\rightarrow \infty}\sup\limits_{x\in {\bf{R}}}
|(v-\overline{v}_{\beta_1,\beta_3},u-\overline{u}_{\beta_1,\beta_3},
\theta-\overline{\theta}_{\beta_1,\beta_3})(x,t)|= 0,
\end{eqnarray}
where the shifts $\beta_1$ and $\beta_3$ are determined by {\rm (2.17)}.
\end{Theorem}

In the rests of this section, we show how the diffusion wave $(v^D,u^D,\theta^D)$ of
 the second characteristic field corresponding to $\lambda_2=0$ is constructed, and also show
some properties of the asymptotic state $(V,U,\Theta)$.
Note that our way to construct the
diffusion wave is much simpler than that in [2]. Keeping in the mind that
$(v^D,u^D,\theta^D)$ is close to $(v_m,u_m,\theta_m)$ and the pressure
$p$ is a Riemann invariant along the second characteristic field, we first
expect that $v^D$ and $\theta^D$ are given in the form
\begin{equation}
v^D = \frac{R}{p_m}\tilde{\Theta},\quad \theta^D = \tilde{\Theta}
\end{equation}
for a function $\tilde{\Theta}$ with the far field states
$\tilde{\Theta}(\pm \infty, t)= \theta_m$,
so that it holds $p^D = \frac{R\theta^D}{v^D} = p_m$.
The the first equation of (1.1) reads as
\begin{equation}
(\frac{R}{p_m}\tilde{\Theta})_t + u^D_x =0,
\end{equation}
and the third equation, after usage of the second one, reads as
\begin{equation}
\frac{R}{\gamma-1}\tilde{\Theta}_t+p_mu^D_x=\frac{\kappa p_m}{R}
\bigg(\frac{\tilde{\Theta}_x}{\tilde{\Theta}}\bigg)_x.
\end{equation}
Substituting (2.28) into (2.29), we have a nonlinear diffusion equation for $\tilde{\Theta}$
\begin{equation}\label{2.14}
\tilde{\Theta}_t= \frac{\kappa (\gamma-1)p_m}{\gamma R^2}\bigg(\frac{\tilde{\Theta}_x}{\tilde{\Theta}}\bigg)_x.
\end{equation}
Since $\Theta^D$ is close to the constant state $\theta_m$, we further approximate
the equation (2.30) by the linear heat equation
\begin{equation}\label{2.14}
\tilde{\Theta}_t= a\,\tilde{\Theta}_{xx},\quad
a= \frac{\kappa (\gamma-1)p_m}{\gamma R^2\theta_m}
= \frac{\kappa (\gamma-1)}{\gamma R v_m},
\end{equation}
and we also impose $\tilde{\Theta}- \theta_m$ carries the mass $\beta_2$, that is,
\begin{equation}
\int (\tilde{\Theta}(x,t)-\theta_m)\,dx = \beta_2.
\end{equation}
Thus, we define $\tilde{\Theta}$ which satisfies (2.31) and (2.32) by employing a typical heat kernel as
\begin{equation}
\tilde{\Theta}(x,t)= \theta_m + \frac{\beta_2}{\sqrt{4\pi
a(1+t)}}e^{-\frac{x^2}{4a(1+t)}},
\end{equation}
and define $u^D$ by the relation (2.28) as
\begin{equation}
u^D = u_m +
\frac{aR}{p_m}{\tilde{\Theta}}_x.
\end{equation}
Then it is straightforward to check that  $(v^D,u^D,\theta^D)$ approximately satisfies (1.1) as in the form
\begin{equation}\label{2.17}
\left\{
\begin{array}{ll}
v^D_t-u^D_x=0,\\[2mm]
u^D_t+p^D_x= \big(\frac{Ra^2}{p_m}{\tilde{\Theta}}_{xx}\big)_x,\\[2mm]
\frac{R}{\gamma-1}\theta^D_t+ (p^Du^D)_x=\kappa\bigg(\dfrac{\theta^D_x}{v^D}\bigg)_x
+ \frac{\kappa p_m}{R}\big(\frac{({\tilde{\Theta}}-\theta_m)}{{\tilde{\Theta}}\theta_m}{\tilde{\Theta}}_x \big)_x.\\[4mm]
(v^D,u^D,\theta^D)(\pm\infty,t)=(v_m,u_m,\theta_m).
\end{array}
\right.
\end{equation}
However, we can not have a divergence form for $E^D$ from (2.35), and so
the condition (2.20) does not hold. In order to overcome this difficulty,
we technically replace the definition of $\theta^D$ by a new one as
\begin{equation}
\theta^D = \tilde{\Theta} - \frac{\gamma-1}{2R}(u^D-u_m)^2
\end{equation}
while the definitions of $v^D$ and $u^D$ are the same as before, that is,
\begin{equation*}
v^D = \frac{R}{p_m}\tilde{\Theta},
\quad u^D = u_m +
\frac{aR}{p_m}{\tilde{\Theta}}_x.
\end{equation*}
Then it follows from (2.36) that
\begin{equation*}
E^D -E_m = (\tilde{\Theta}-\theta_m) +\frac{\gamma-1}{R}(u^D-u_m)u_m
\end{equation*}
which implies
\begin{equation}
\int (E^D(x,t) -E_m) \,dx = \beta_2.
\end{equation}
Then, direct calculations show that the new $(v^D,u^D,\theta^D)$ satisfies
\begin{equation}\label{2.17}
\left\{
\begin{array}{ll}
v^D_t-u^D_x=0,\\[2mm]
u^D_t+p^D_x= (R_1^D)_x,\\[2mm]
\bigg(\frac{R}{\gamma-1}\theta^D+\frac{(u^D)^2}{2}\bigg)_t+ (p^Du^D)_x=\kappa\bigg(\dfrac{\theta^D_x}{v^D}\bigg)_x
+ (R_2^D)_x, \\[4mm]
(v^D,u^D,\theta^D)(\pm\infty,t)=(v_m,u_m,\theta_m),\\[4mm]
\int (m^D - m_m) \,dx = \beta_2 r_2,
\end{array}
\right.
\end{equation}
where the remainder terms $R_1^D$ and $R_2^D$ are given by
\begin{equation}
R_1^D = \frac{Ra^2}{p_m}{\tilde{\Theta}}_{xx}
-\dfrac{\gamma-1}{2Rv^D}(u^D-u_m)^2,
\end{equation}
and
\begin{equation}
R_2^D = \frac{\kappa p_m}{R}\big(\frac{{\tilde{\Theta}}-\theta_m}
{{\tilde{\Theta}}\theta_m}\big){\tilde{\Theta}}_x
-\dfrac{\gamma-1}{2v^D}(u^D-u_m)^2U
+ \frac{Ra^2u_m}{p_m}{\tilde{\Theta}}_{xx}
+ \dfrac{\kappa(\gamma-1)}{Rv^D}(u^D-u_m)(u^D-u_m)_x.
\end{equation}
By the definition of $(v^D,u^D,\theta^D)$, it is also noted that
$R_1^D$ and $R_2^D$ satisfy the following pointwise estimates as $\beta_2 \to 0$:
\begin{equation}
|(\frac{\partial }{\partial x})^i(R^D_1)|\le
\frac{C\beta_2}{(1+t)^{2}}e^{-\frac{x^2}{4a(1+t)}},\quad
|(\frac{\partial }{\partial x})^i(R^D_2)|\le
\frac{C\beta_2}{(1+t)^{\frac{3}{2}}}e^{-\frac{x^2}{4a(1+t)}},\quad
\ x \in {\bf R},\ t\ge 0,\ 0\le i \le 3.
\end{equation}

Finally, let us see several properties of the asymptotic state $Z(x,t)={}^t(V,U,\Theta)(x,t)$
defied by (2.22). By using the equations of the viscous shock waves and diffusion wave,
we firstly can see $Z(x,t)={}^t(V,U,\Theta)(x,t)$ approximately satisfies (1.1) as
in the form
\begin{eqnarray}\label{2.42}
\left\{
\begin{array}{ll}
V_t-U_x=0,\\[2mm]
U_t+P_x= (R_{1})_x,\\[2mm]
\bigg(\frac{R}{\gamma-1}\Theta+\frac{U^2}{2}\bigg)_t
+(PU)_x=\bigg(\dfrac{\kappa \Theta_x}{V}\bigg)_x+ (R_{2})_x,\\[4mm]
(V,U,\Theta)(\pm\infty,t)=(v_\pm,u_\pm,\theta_\pm),
\end{array}
\right.
\end{eqnarray}
where $R_1$, $R_2$ are given by
\begin{eqnarray}
R_{1}&=& P- (P_1+P_3-p_m +(p^D-p_m)) + R^D_1,\\[2mm]
R_{2}&=&-\frac{\kappa \Theta_x}{V}+
\frac{\kappa \Theta_{1,x}}{V_1}
+\frac{\kappa  \Theta_{3,x}}{V_3}+ \frac{\kappa \theta^D_x}{v^D} \nonumber\\[2mm]
&& + PU-(P_1U_1+P_3U_3-p_mu_m + (p^Du^D-p_mu_m)) + R^D_2.\nonumber
\end{eqnarray}
Noting that
$$
p^D-p_m = -\frac{\gamma -1}{2v^D}(u^D-u_m)^2
$$
and
\begin{eqnarray*}
\Theta &=& \Theta_1 + \Theta_3 -\theta_m + (\theta^D-\theta_m)\\
&-& \frac{\gamma -1}{R}\big((U_1-u_m)(U_3-u_m)+(U_1-u_m)(u^D-u_m)+(U_3-u_m)(u^D-u_m) \big),
\end{eqnarray*}
we can show the estimate
\begin{equation}
|R_1| \leq C(|u^D-u_m|^2+ |Z_1-z_m||Z_3-z_m|+ (|Z_1-z_m|+|Z_3-z_m|)|z^D-z_m|) + |R^D_1|.
\end{equation}
The wave interaction terms in (2.44) for small $\delta$ and $\beta_2$
can be estimated  by Proposition 2.1 and the definition of
the diffusion wave as in the same way as in [2]:
\begin{eqnarray}
|Z_1-z_m||Z_3-z_m|&\leq & C\delta_1\delta_3(e^{-c\delta_1(|x|+t)+c\delta_1|\beta_1|}+e^{-c\delta_3(|x|+t)+c\delta_3|\beta_3|}) \nonumber \\[2mm]
&\leq & C\delta^2e^{-c\delta(|x|+t)},
\end{eqnarray}
and for $i=1,3$
\begin{eqnarray}\label{2.22}
|Z_i-z_m||\tilde{\Theta}|&\leq& C\delta_ie^{-c\delta_i(|x|+t)+c\delta_i|\beta_i|}\frac{|\beta_2|}{(1+t)^\frac{1}{2}}e^{-\frac{cx^2}{1+t}}
+C\delta_i|\beta_2|e^{-c(|x|+t)}\nonumber\\[2mm]
&\le& C|\beta_2|\delta^{\frac32}e^{-c\delta(|x|+t)}
+ C\frac{|\alpha_2|}{(1+t)^{\frac32}}e^{-\frac{cx^2}{1+t}}
+C(\delta+|\beta_2|)e^{-c(|x|+t)},
\end{eqnarray}
where we used the fact that $\delta_1\beta_1$ and $\delta_3\beta_3$ are uniformly bounded by (\ref{2.13})
for small $\delta_1, \delta_3$ as long as the initial perturbation stay bounded.
By (2.41), (2.46) and (2.47), we have
\begin{equation}
|R_1| \leq
C(\delta^2+|\beta_2|\delta^{\frac32})e^{-c\delta(|x|+t)}+
C\frac{|\beta_2|}{(1+t)^{\frac32}}e^{-\frac{cx^2}{1+t}}+C(\delta+|\beta_2|)e^{-c(|x|+t)}.
\end{equation}
In what follows, let us denote $A \approx B$ when it holds
$$
|A-B|\le
 \
C(\delta^2+|\beta_2|\delta^{\frac32})e^{-c\delta(|x|+t)}+
C\frac{|\beta_2|}{(1+t)^{\frac32}}e^{-\frac{cx^2}{1+t}}+C(\delta+|\beta_2|)e^{-c(|x|+t)}.
$$
Then, we can simply rewrite (2.47) as $R_1 \approx 0$.
Similarly, we can show that
\begin{equation}
(\frac{\partial }{\partial x})^i(R_1),\
(\frac{\partial }{\partial x})^i(R_2) \approx 0, \quad
0\le i \le 3.
\end{equation}

\section{Reformulation of problem and local existence}
\setcounter{equation}{0}
As in the previous paper [2], we first reformulate the Cauchy problem
(1.1),(1.3) around the approximate solution $(V,U,\Theta)$
given by (2.22) to one for an integral system of (1.1).
Introducing the anti-derivative variables which correspond to the density,
velocity, total energy and absolute temperature as
\begin{eqnarray}\label{3.1}
\Phi(x,t)&=& \int_{-\infty}^x{(v-V)(y,t)}\,dy,\nonumber\\[2mm]
\Psi(x,t)&=& \int_{-\infty}^x{(u-U)(y,t)}\,dy,\\[2mm]
\overline{W}(x,t)&=& \int_{-\infty}^x{\left((\frac{R}{\gamma-1}\theta+\frac{u^2}{2})-
(\frac{R}{\gamma-1}\Theta+\frac{U^2}{2})\right)(y,t)}\,dy,\nonumber\\[2mm]
W(x,t)&=& \frac{\gamma-1}{R}(\overline{W}-U\Psi)(x,t), \nonumber
\end{eqnarray}
we look for the solution in the form
\begin{equation}
v-V=\Phi_x,\quad u-U=\Psi_x,\quad
\theta-\Theta=W_x+\frac{\gamma-1}{R}(U_x\Psi-\frac{1}{2}\Psi^2_x),
\end{equation}
where we expect $(\Phi,\Psi,W) \in C([0,\infty);H^3)$.
Then, substituting (3.2) to (1.1), subtracting (\ref{2.42}), and integrating the resulting system
with respect to $x$, we have the
integrated system for $(\Phi,\Psi,W)$
\begin{eqnarray}\label{3.3}
\left\{
\begin{array}{ll}
\Phi_t-\Psi_x=0, \\[2mm]
\Psi_t+R\left(\dfrac{\Theta+W_x+\frac{\gamma-1}{R}(U_x\Psi-\frac{1}{2}\Psi^2_x)}{V+\Phi_x}
-\dfrac{\Theta}{V}\right)=-R_1, \\[5mm]
\frac{R}{\gamma-1}W_t+U_t \Psi+
R\left(\dfrac{\Theta+W_x+\frac{\gamma-1}{R}(U_x\Psi-\frac{1}{2}\Psi^2_x)}{V+\Phi_x}\right)\Psi_x\\[5mm]
\qquad = \kappa \left(\dfrac{\big(\Theta+W_x+\frac{\gamma-1}{R}(U_x\Psi-\frac{1}{2}\Psi^2_x)\big)_x}{V+\Phi_x}
-\dfrac{\Theta_x}{V}\right)
-R_{2}+UR_1,
\end{array}
\right.
\end{eqnarray}
and the initial data should be given by
\begin{equation}
(\Phi,\Psi,W)(0)=(\Phi_0,\Psi_0,W_0) \in H^3
\end{equation}
where $W_0(x) = \frac{\gamma-1}{R}(\overline{W}_0(x)-U(x,0)\Psi_0(x))$,
$(\Phi_0,\Psi_0,\overline{W}_0)$ are given by (2.25).
Now we focus our attention on the Cauchy problem to the reformulated system (\ref{3.3})
with initial date (3.4), and construct the global solution in time
$(\Phi,\Psi,W)\in C([0,\infty);H^3)$ under some smallness conditions
on $(\Phi_0,\Psi_0,W_0)$ and $\delta+|\beta_2|$,
combining the local existence in time and the {\it a priori} estimates.
To consider the local existence of the solution, since
the system (3.3) is non-autonomous, we need to consider the Cauchy problem for
(3.3) where the initial date is given at general time $\tau\ge 0$ as
$$
(\Phi,\Psi,W)(\tau)=(\Phi_\tau,\Psi_\tau,W_\tau) \in H^3.
\eqno{(3.4)_\tau}
$$
In order to state the local existence precisely,
we define the solution set $X(I)$ for any interval $I\subseteq \bf{R}_+$ by
$$
X(I)=\{ (\Phi,\Psi,W)\in C(I;H^3)|\ (\Phi_x,\Psi_x, W_x)\in L^2(I;H^2),
\xi \in L^2(I;H^3)\},
$$
where $\xi = W_x+\frac{\gamma-1}{R}(U_x\Psi-\frac{1}{2}\Psi^2_x)$.
We further choose a positive constant $\bar{\delta}_0$ for given $(v_-,u_-,\theta_-)$ such that
$\delta+|\beta_2|\le \bar{\delta}_0$ implies
\begin{equation}
\sup_{x\in {\bf R}, t\ge 0}|(V-v_-,U-u_-,\Theta-\theta_-)(x,t)| \le \frac{1}{4}\min (v_-,\theta_-),
\end{equation}
and also choose a positive constant $\bar{\epsilon}_0\,(\le \frac{1}{2}\min (v_-,\theta_-))$ such that
$\|(\Phi,\Psi,W)(t)\|_2 \le \bar{\epsilon}_0$ implies
\begin{equation}
\sup_{x\in {\bf R}}|(\Phi,\Psi,W_x+\frac{\gamma-1}{R}(U_x\Psi-\frac{1}{2}\Psi^2_x))(x,t)|
\le \frac{1}{2}\min (v_-,\theta_-).
\end{equation}

\begin{Proposition}{\rm (Local existence)} \
For any fixed $(v_-,u_-,\theta_-)$, there exist positive constants
$\bar{\epsilon}_1\,(\le \bar{\epsilon}_0)$
and $\bar{C}_1\,(\bar{C}_1\bar{\epsilon}_1\le \bar{\epsilon}_0)$ such that
the following statements hold. Under the assumption
$\delta+|\beta_2|\le \bar{\delta}_0$,
for any constant $M\in (0,\bar{\epsilon}_1)$,
there exists a positive constant
$t_0=t_0(M)$ not depending on $\tau$
such that if $\|(\Phi_\tau,\Psi_\tau,W_\tau)\|_{H^3}\le M$, then
the Cauchy problem $(3.3),(3.4)_\tau$ admits a unique solution
 $(\Phi,\Psi,W) \in X([\tau, \tau+t_0])$ satisfying
\begin{eqnarray*}
\sup\limits_{t\in [\tau,\tau+t_0]}\|(\Phi,\Psi,W)(t)\|_3 \leq \bar{C}_1 M.
\end{eqnarray*}
\end{Proposition}
Here we should note that since $\bar{C}_1 M\le \bar{C}_1\bar{\epsilon}_1\le \bar{\epsilon}_0$,
the Sobolev's inequality and (3.5),(3.6) imply the local solutions constructed above
satisfy
\begin{equation}
\frac{1}{4}v_-\le (V+\Phi_x)(x,t) \le \frac{7}{4}v_-,\quad
\frac{1}{4}\theta_-\le (\Theta+\xi)(x,t) \le \frac{7}{4}\theta_-,
\end{equation}
where $\xi = W_x+\frac{\gamma-1}{R}(U_x\Psi-\frac{1}{2}\Psi^2_x)$.
The proof of Proposition 3.1 is not so trivial because the system (3.3) is
purely nonlinear, despite in the previous viscous and heat-conductive case
in [2] the corresponding integral system is quasi-linear.
Therefore we  will simply show a strategy to prove the local
existence at the last of this section.
In order to have the global solution by using the
local existence, Proposition 3.1, repeatedly, we need to establish
the {\it a priori} estimates for the solution of (3.3),(3.4).
Setting
$$
N(T)=\sup\limits_{t\in [0,T]}\|(\Phi,\Psi,W)(t)\|_3,
$$
we shall show the following {\it a priori} estimates in
the next section.
\begin{Proposition}{\rm({\it a priori} estimates)} \
Under the same assumptions in {\rm Theorem 2.1},
there exist positive constants
$\delta_0\,(\le \bar{\delta}_0)$,
$\epsilon_0\,(\le \bar{\epsilon}_1)$
and $C_0$ which depend only on $(v_-,u_-.\theta_-)$ such that
the following statements hold. If $(\Phi,\Psi,W)\in X([0,T])$ is the solution of
$(3.3),(3.4)$ for some $T>0$, $\delta+|\beta_2|\leq \delta_0 $, and
$N(T)\leq \epsilon_0$, then it holds for $t\in [0,T]$,
\begin{eqnarray}
&&\|(\Phi,\Psi,W)(t)\|^2_3+\int^t_0\|(|U_{1x}|+|U_{3x}|)^{\frac 1 2}(\Psi,W)(\tau)\|^2\,
d\tau\\[2mm]
&&+\int^t_0{\big(\|(\Phi_x,\Psi_x,W_x)(\tau)\|^2_2+\|\xi(\tau)\|^2_3\big)}\,d\tau
\leq C_0(\|(\Phi_0,\Psi_0,W_0)\|^2_3+\delta^{\frac{1}{2}}+|\beta_2|),\nonumber
\end{eqnarray}
where $\xi = W_x+\frac{\gamma-1}{R}(U_x\Psi-\frac{1}{2}\Psi^2_x)$.
\end{Proposition}
Once  Proposition 3.2 is proved, choosing $\delta$, $|\beta_2|$ and
$\|(\Phi_0,\Psi_0,W_0)\|^2_3$ suitably small such as
$$
\delta+|\beta_2|\leq \delta_0,\quad
\|(\Phi_0,\Psi_0,W_0)\|_3\le \frac{\epsilon_0}{\bar{C}_1},\quad
C_0(\|(\Phi_0,\Psi_0,W_0)\|^2_3+\delta^{\frac{1}{2}}+|\beta_2|)\le
\left(\frac{\epsilon_0}{\bar{C}_1}\right)^2,
$$
we can construct the global solution $(\Phi,\Psi,W)\in X([0,\infty))$
by combining Propositions 3.1 and 3.2 as in the previous papers, and
can show the estimate (3.8) holds for all $t\ge 0$.
Furthermore, we can see from the estimate  (3.8) and the system (\ref{3.3}) that
\begin{eqnarray*}
\int^\infty_0{\|(\Phi_x,\Psi_x,W_x)(t)\|^2}dt+
\int^\infty_0{\bigg|\frac{d}{dt}\|(\Phi_x,\Psi_x,W_x)(t)\|^2\bigg|}dt<+\infty,
\end{eqnarray*}
which together with the Sobolev's inequality leads to the asymptotic behavior of the solution
\begin{eqnarray*}
\lim\limits_{t\rightarrow \infty}\sup\limits_{x\in {\bf R}}|(\Phi,\Psi,W,\Phi_x,\Psi_x,W_x)(x,t)|=0.
\end{eqnarray*}
Noticing the original unknown variables $(v,u,\theta)$ can be restored by
$$
v= V + \Phi_x,\quad u= U+\Psi_x, \quad \theta=\Theta+W_x+\frac{\gamma-1}{R}
(U_x\Psi-\frac{1}{2}\Psi^2_x),
$$
and also the smallness of $\delta$, $|\beta_2|$ and
$\|(\Phi_0,\Psi_0,W_0)\|^2_3$ are controlled by the
smallness of $\delta$ and $I(v_0,u_0,\theta_0)$ under the
assumptions of Theorem 2.1, we can have the global solution
in Theorem 2.1 with the asymptotic behavior
\begin{equation}
\lim\limits_{t\rightarrow \infty}\sup\limits_{x\in {\bf{R}}}
|(v-V,u-U,
\theta-\Theta)(x,t)|= 0.
\end{equation}
Since the diffusion wave $z^D$ uniformly tends toward the
constant state $z_m$ in the definition of $(V,U,\Theta)$ in (2.22),
we finally have the desired asymptotic behavior
$$
\lim\limits_{t\rightarrow \infty}\sup\limits_{x\in {\bf{R}}}
|(v-\overline{v}_{\beta_1,\beta_3},u-\overline{u}_{\beta_1,\beta_3},
\theta-\overline{\theta}_{\beta_1,\beta_3})(x,t)|= 0.
$$
Thus Theorem 2.1 can be proved by Propositions 3.1 and 3.2.

Now we turn to Proposition 3.1 and show a strategy of the proof.
We treat only the case $\tau=0$ for simplicity because the
case $\tau>0$ is just similar. Since the system (3.3) is fully nonlinear,
differentiating (3.3) with respect to $x$, and introducing
the new variables $(\phi,\psi,\xi)$ by
\begin{eqnarray}\label{3.9}
&&\phi=\Phi_x,\quad \psi=\Psi_x, \quad
 \xi=W_x+\frac{\gamma-1}{R}(U_x\Psi-\frac{1}{2}\Psi^2_x),
\end{eqnarray}
we have the system in terms of $(\phi,\psi,\xi)$
\begin{eqnarray}\label{3.10}
\left\{
\begin{array}{ll}
\phi_t-\psi_x=0, \\[2mm]
\psi_t+\left(\frac{R\xi}{V+\phi}-\frac{P\phi}{V+\phi}\right)_x=-R_{1,x},\\[4mm]
\frac{R}{\gamma-1}\xi_t+\frac{R(\Theta+\xi)}{V+\phi}\psi_x
+\left(\frac{R\xi}{V+\phi}-\frac{P\phi}{V+\phi}\right)U_x=
\kappa\left(\frac{\xi_x}{V+\phi}-\frac{\Theta_x\phi}{V(V+\phi)}\right)_x
+U R_{1,x}-R_{2,x},
\end{array}
\right.
\end{eqnarray}
with the initial data
\begin{equation}
(\phi,\psi,\xi)(0)= (\phi_0,\psi_0,\xi_0):=
\left(\Phi_{0x},\Psi_{0x}, W_{0x}
+\frac{\gamma-1}{R}(U_x(\cdot,0)\Psi_0-\frac{1}{2}\Psi^2_{0x})\right) \in H^2.
\end{equation}
We can also rewrite the system (3.3)
as in the form
\begin{eqnarray}
\left\{
\begin{array}{ll}
\Phi_t=\psi,\\[2mm]
\Psi_t=-\left(\frac{R\xi}{V+\phi}-\frac{P\phi}{V+\phi}\right)-R_1, \\[3mm]
\frac{R}{\gamma-1}W_t+U_t \Psi=
\kappa\left(\frac{\xi_x}{V+\phi}-\frac{\Theta_x\phi}{V(V+\phi)}\right)
-\frac{R(\Theta+\xi)}{V+\phi}\psi+UR_1-R_{2},
\end{array}
\right.
\end{eqnarray}
with the initial data
\begin{equation}
(\Phi,\Psi,W)(0)=(\Phi_0,\Psi_0, W_0) \in H^3.
\end{equation}
Then, we consider the Cauchy problem for the extended system of (3.11),(3.13)
for $(\Phi,\Psi,W,\phi,\psi,\xi)$ with the initial data (3.12),(3.14).
Keeping in mind the system (3.11) is closed in terms of $(\phi,\psi,\xi)$, we
rewrite it again as in a form of quasi-linear hyperbolic-parabolic system
\begin{eqnarray}
\left\{
\begin{array}{ll}
\phi_t-\psi_x=0, \\[3mm]
\psi_t-\frac{R(\Theta+\xi)}{(V+\phi)^2}\phi_x=F_1(\phi,\psi,\xi_x),
\quad \quad x\in {\bf R}, t\ge 0,\\[3mm]
\frac{R}{\gamma-1}\xi_t-
\frac{\kappa }{V+\phi}\xi_{xx}=F_2(\phi,\xi,\phi_x, \psi_x,\xi_x),
\end{array}
\right.
\end{eqnarray}
where
\begin{eqnarray}\label{3.12}
\begin{array}{ll}
F_1(\phi,\psi,\xi_x)=\left(\frac{P_x}{V+\phi}
-\frac{PV_x}{(V+\phi)^2}\right)\phi
+\frac{RV_x \xi}{(V+\phi)^2}-\frac{R \xi_x}{V+\phi}-R_{1,x},\\[3mm]
F_2(\phi,\psi,\xi_x)=-\frac{R(\Theta+\xi)}{V+\phi}\psi_x
+\frac{P\phi-R\xi}{V+\phi}U_x\\[3mm]
\qquad \qquad \quad\quad-\frac{\kappa\xi_x(V_x+\phi_x)}{(V+\phi)^2}
-\kappa\left(\frac{\Theta_x\phi}{V(V+\phi)}\right)_x+UR_{1,x}-R_{2,x}.
\end{array}
\end{eqnarray}
We note that as long as we look for the local solution of (3.15),(3.12) for $t\in [0,t_0]$
satisfying
$$
\sup_{t\in [0,t_0]}\|(\phi,\psi,\xi)(t)\|_2 \le \bar{\epsilon}_0,
$$
the principal part of the left hand side of (3.15)
for $(\phi,\psi)$ is strictly hyperbolic and that for $\xi$ is uniformly parabolic because of
(3.7), and the left hand side of (3.15) can be regarded as lower order terms.
Therefore, applying the arguments on such hyperbolic-parabolic systems by Kawashima [6],[7],
we can show the existence and uniqueness of the local solution $(\phi,\psi,\xi)\in
C([0,t_0];\,H^2)$ satisfying $\xi \in  L^2(0,t_0;\, H^3)$ for a suitably small
$t_0=t_0(M)>0$, provided $\|(\Phi_0,\Psi_0,W_0)\|_3 \le M$ for small $M$,
say $M \le \bar{\epsilon}_1$. Then, plugging $(\phi,\psi,\xi)$ into (3.13), we can easily
have the solution of (3.13),(3.14) satisfying $(\Phi,\Psi) \in C([0,t_0];\,H^2)$ and
$W \in C([0,t_0];\,H^1)\cap L^2(0,t_0;\,H^2)$. Now, we can show the relations (3.10) hold,
that is,
\begin{equation}
\phi=\Phi_x,\quad \psi=\Psi_x, \quad
 \xi=W_x+\frac{\gamma-1}{R}(U_x\Psi-\frac{1}{2}\Psi^2_x).
\end{equation}
In fact, by (3.11)-(3.14), we can see
\begin{eqnarray}
\left\{
\begin{array}{ll}
(\phi-\Phi_x)_t = 0, \\[2mm]
(\psi-\Psi_x)_t = 0,\\[2mm]
\left((\frac{R}{\gamma-1}\xi+U \psi+\frac{1}{2}\psi^2)
-(\frac{R}{\gamma-1}W+U \Psi)_x\right)_t=0.
\end{array}
\right.
\end{eqnarray}
and the initial data satisfy
\begin{equation}
(\phi-\Phi_x)(0)=(\psi-\Psi_x)(0)=
\big((\frac{R}{\gamma-1}\xi+U \psi+\frac{1}{2}\psi^2)
-(\frac{R}{\gamma-1}W+U \Psi)_x\big)(0)=0.
\end{equation}
Hence, it follows from (3.18) and (3.19) that
\begin{equation}
\phi=\Phi_x,\quad \psi=\Psi_x, \quad
\frac{R}{\gamma-1}\xi+U \psi+\frac{1}{2}\psi^2
=(\frac{R}{\gamma-1}W+U \Psi)_x,
\end{equation}
which is nothing but (3.17). By (3.17), we eventually can conclude
$(\Phi,\Psi,W) \in C([0,t_0];\,H^3)$ and
$(\Phi,\Psi,W)$ gives the desired local solution of
$(3.3),(3.4)_0$. Thus, we can show Proposition 3.1.

\section{ {\it A priori} estimates}
\setcounter{equation}{0}

In this section, we show the {\it a prior} estimates, Proposition 3.2.
Although the proof is basically given by following the procedure as in
[2], we need more subtle estimates than [2] because of lack of
viscosity term. Throughout of this section, we assume that
 $(\Phi,\Psi,W)\in X([0,T])$ is a solution of (\ref{3.3}) for some $T>0$,
$\delta+|\beta_2| \le \delta_0\,(\le \bar{\delta}_0)$, and
$N(T) \le \epsilon_0\,(\le \bar{\epsilon}_1)$, where $\epsilon_0$ and $\delta_0$ will
be chosen suitably small later.  Then it is noted that as in (3.7)
$V,V+\Phi_x,\Theta,\Theta+\xi$  are
uniformly positive  on $[0,T]$ by Sobolev's inequality as
$$
\inf_{x,t}V\ge \frac{3v_-}{4},\quad
\ \inf_{x,t}(V+\Phi_x) \ge \frac{v_-}{4},\quad
\inf_{x,t}\Theta \ge \frac{3\theta_-}{4},\quad
\ \inf_{x,t}(\Theta+\xi)\ge \frac{\theta_-}{4}.
$$
For the {\it a prior} estimates, we next rewrite the system (3.3) again
so that all linearized parts at $(\Phi,\Psi,W)$ are collected in left hand
side, and nonlinear and inhomogeneous terms in right hand side as
\begin{eqnarray}
\left\{
\begin{array}{ll}
\Phi_t-\Psi_x=0, \\[2mm]
\Psi_t-\frac{P}{V}\Phi_x+\frac{R}{V}W_x
+\frac{\gamma-1}{V}U_x \Psi=N_1-R_1, \\[2mm]
\frac{R}{\gamma-1}W_t+P \Psi_x+U_t \Psi-
\frac{(\gamma-1)\kappa}{RV}(U_x\Psi)_x+
\frac{\kappa \Theta_x}{V^2}\Phi_x-
\frac{\kappa}{V}W_{xx}=N_2-R_2+UR_1,
\end{array}
\right.
\end{eqnarray}
where
\begin{eqnarray*}
N_1&=&\frac{\gamma-1}{2V}\Psi^2_x-\bigg(p-P
+\frac{P}{V}\Phi_x-\frac{R}{V}(\theta-\Theta) \bigg),\\[2mm]
N_2&=&(P-p)\Psi_x-\frac{(\gamma-1)}{RV}\Psi_x\Psi_{xx}-
\frac{\kappa }{vV}(\theta-\Theta)_x\Phi_x
+\frac{\kappa }{vV^2}\Theta_x\Phi_x^2,
\\[2mm]
p&=&\frac{R\theta}{v},\quad v= V+\Phi_x,\quad
\theta=\Theta+W_x+\frac{\gamma-1}{R}(U_x\Psi-\frac{1}{2}\Psi^2_x).
\end{eqnarray*}
It is noted that for small $\|(\Phi,\Psi,W)\|_3$
\begin{eqnarray}
N_1&=&O(1)(\Phi^2_x+\Psi^2_x+W^2_x+|U_x|\Psi^2),\\[2mm]
N_2&=&O(1)(\Phi^2_x+\Psi^2_x+W^2_x+|U_x|\Psi^2+|\Psi_x||\Psi_{xx}|+|\Phi_x||(\theta-\Theta)_x|).
\nonumber
\end{eqnarray}
Then the proof of the Proposition 3.2 is given by
the following series of lemmas.

\begin{Lemma} If $\delta_0$ and $\epsilon_0$ are suitably small, it holds that
\begin{eqnarray}\label{4.2}
&&\|(\Phi,\Psi,W)(t)\|^2+\int^t_0\int
\big((|U_{1x}|+|U_{3x}|)(\Psi^2+W^2)+W_{x}^2\big)\,dxd\tau
\\[2mm]
&\leq &C(\|(\Phi_0,\Psi_0,W_0)\|^2+\delta^\frac{1}{2}+|\beta_2|)
+C(\epsilon_0+\delta_0)\int^t_0\int{(\Phi_{x}^2
+\Psi^2_x+\xi^2_x+\Psi^2_{xx})}\,dxd\tau,\nonumber
\end{eqnarray}
where $\xi=W_x+\frac{\gamma-1}{R}(U_x\Psi-\frac{1}{2}\Psi^2_x)$.
\end{Lemma}
{\bf Proof.}\ \ Multiplying the equation $(4.1)_1$ by $\Phi$, $(4.1)_2$ by $\frac{V}{P}\Psi$ and
$(4.1)_3$ by $ \frac{R}{P^2}W$ respectively,  and adding all the resultant equations, we
obtain
\begin{eqnarray}
&&\bigg\{\frac{\Phi^2}{2}+\frac{V}{2P}\Psi^2+\frac{R^2}{2(\gamma-1)P^2}W^2\bigg\}_t
+\bigg\{-\Phi\Psi+\frac{R}{P}\Psi W-\frac{(\gamma-1)\kappa}{VP}U_x\Psi W-
\frac{R\kappa}{VP^2}WW_x\bigg\}_x\nonumber\\[2mm]
&&+A\Psi^2+J_1+J_2=\frac{V}{P}\Psi (N_1- R_1)+ \frac{R}{P^2}W (N_2- R_2+UR_1).
\end{eqnarray}
where
\begin{eqnarray*}
A&=&-\big(\frac{V}{2P}\big)_t+\frac{\gamma-1}{P}U_x,\nonumber\\[2mm]
J_1&=&\frac{R\kappa}{VP^2}W^2_x+
\big(\frac{R\kappa}{VP^2}\big)_xWW_x+\frac{R^2P_t}{(\gamma-1)P^3}W^2,\\[2mm]
J_2&=&\big(\frac{(\gamma-1)\kappa}{VP^2}W\big)_xU_x\Psi
+\frac{\kappa \Theta_x}{V^2}\frac{R}{P^2}W\Phi_x-\frac{R}{P^2}R_{1x}W\Psi.
\end{eqnarray*}
We estimate the terms $A,J_1$ and $J_2$ one by one.
Since it holds
$$
P_t\approx P_{1t}+P_{3t}, \quad P_x\approx P_{1x}+P_{3x},\quad U_x\approx U_{1x}+U_{3x},
$$
we obtain
\begin{eqnarray*}
A\approx \big(-\big(\frac{V_1}{2P_1}\big)_t+\frac{\gamma-1}{P_1}U_{1x}\big)
+\big(-\big(\frac{V_3}{2P_3}\big)_t+
\frac{\gamma-1}{P_3}U_{3x}\big)=:A_1+A_3.
\end{eqnarray*}
Using the fact $P_{it}=-s_iP_{ix}=s_iU_{it}=s_i^2(-U_{ix})>0$,
we further obtain from  Proposition 2.1 that
\begin{eqnarray*}
A_i&=&-\frac{V_{it}}{2P_i}+\frac{V_iP_{it}}{2P^2_i}+\frac{\gamma-1}{P_i}U_{ix}
=\frac{(-U_{ix})}{2P^2_i}\big( (s^2_iV_i-\gamma P_i)+(3-\gamma)P_i\big)\\[2mm]
&\geq& c|U_{ix}|((3-\gamma)p_m-C\delta_i).
\end{eqnarray*}
Since $\gamma \in [1,3)$, choosing $\delta_0$ suitably small, we get
\begin{eqnarray}\label{4.4}
A\Psi^2 \geq c(|U_{1x}|+|U_{3x}|)\Psi^2 - \tilde{R}\Psi^2
\end{eqnarray}
for a positive constant $c$,
where and also in what follows $\tilde{R}$ represents some error functions which
satisfy $\tilde{R}\approx 0$.
For $J_1$, we first estimate the second term in $J_1$ by the Cauchy inequality as
\begin{eqnarray*}
|\big(\frac{R\kappa}{VP^2}\big)_xWW_x| \leq \frac{R\kappa}{2VP^2}W^2_x
+C(|U_{1x}|^2+|U_{3x}|^2+(\theta^D_x)^2)W^2.
\end{eqnarray*}
Then, noting $P_{it}=s_i^2(-U_{ix})>0$ and $(\theta^D_x)^2\approx0$,
we have for a positive constant $c$
\begin{eqnarray}
J_1 \geq c(W^2_x +(|U_{1x}|+|U_{3x}|)W^2) - \tilde{R}W^2.
\end{eqnarray}
For $J_2$, using $(\theta^D_x)^2\approx0$ and $R_{1,x}\approx 0$, we similarly can get
\begin{eqnarray}
J_2 \geq -C\delta_0( W^2_x+\Phi^2_x) -C\delta_0(|U_{1x}|+|U_{3x}|)(\Psi^2+W^2)
-\tilde{R}(\Psi^2+W^2).
\end{eqnarray}
Next we estimate the terms $\Psi N_1$ and $ W N_2$. From (4.2) and Sovolev's inequality,
it easily follows
\begin{eqnarray}
|\Psi N_1|&\leq& C\epsilon_0 (\Phi^2_x+\Psi^2_x+W^2_x+(|U_{1x}|+|U_{3x}|)\Psi^2)
+\tilde{R}\Psi^2,\\[2mm]
|W N_2|&\leq& C\epsilon_0(\Phi^2_x+\Psi^2_x+W^2_x+\xi^2_x+\Psi^2_{xx}
+(|U_{1x}|+|U_{3x}|)\Psi^2)+\tilde{R}\Psi^2.\nonumber
\end{eqnarray}
Finally, recalling $R_1\approx 0$ and $R_2\approx 0$ by (2.49),
we deal with all the error terms arising from the relation $``\approx"$
like $\tilde{R}\Psi W_x, \tilde{R} \Psi $
and $\tilde{R}W$ appeared in (4.4) and $(\Psi^2 +W^2)\tilde{R}$ in (4.5)-(4.8).
We can see that
all the integrations of such terms on ${\bf{R}}\times (0,t)$ are estimated by
\begin{equation}
C\delta_0\int^t_0{\|(\Psi_x,W_x)(\tau)\|^2}\,d\tau+C(\delta^{\frac{1}{2}}+|\beta_2|).
\end{equation}
Here we show (4.9) only for typical terms $\tilde{R}\Psi W_x$ and $\tilde{R} \Psi$
because the other terms are similarly estimated:
\begin{eqnarray}
\int^t_0\int{|\tilde{R}|(|\Psi W_x|+|\Psi|)}\,dxd\tau\leq
  C\delta_0\int^t_0{\|W_x\|^2}\,d\tau+C\int^t_0\int{|\tilde{R}||\Psi|}\,dxd\tau,
\end{eqnarray}
and
\begin{eqnarray}
&&\int^t_0\int{|\tilde{R}||\Psi|}\,dxd\tau\leq
C\int^t_0\int{(\delta^2+|\beta_2|\delta^{\frac{3}{2}})e^{-c\delta(|x|+\tau)}|\Psi|}\,dxd\tau
\nonumber\\[2mm]
&&+C\delta\int^t_0\int{\frac{|\beta_2|}{(1+\tau)^\frac{3}{2}}
 e^{-\frac{cx^2}{1+\tau}}|\Psi|}\,dxd\tau
+C\int^t_0\int{(\delta+|\beta_2|)e^{-c(|x|+\tau)}|\Psi|}\,dxd\tau \\[2mm]
 &\leq&C\int^t_0{(\delta^\frac{3}{2}+|\beta_2|\delta)e^{-c\delta\tau}
 \|\Psi\|}d\tau+ C\delta\int^t_0{\frac{C|\beta_2|}{(1+\tau)^\frac{5}{4}}
 \|\Psi\|}d\tau+C(\delta+|\beta_2|) \nonumber \\[2mm]
  &\leq& C(\delta^{\frac{1}{2}}+|\beta_2|)+C(\delta+|\beta_2|)\leq
C(\delta^{\frac{1}{2}}+|\beta_2|).\nonumber
\end{eqnarray}
Combing all the estimates (4.5)-(4.11), integrating (4.4) on ${\bf{R}} \times (0,t)$ and choosing
$\epsilon_0$ and $\delta_0$
suitably small, we have the desired estimate (4.3). This completes the proof of Lemma 4.1.

\medskip

For the estimates for higher derivatives, using the variables
$$
\phi=\Phi_x,\quad \psi=\Psi_x, \quad
 \xi=W_x+\frac{\gamma-1}{R}(U_x\Psi-\frac{1}{2}\Psi^2_x),
$$
we write the equation (3.11) in terms of $(\phi,\psi,\xi)$  as
\begin{eqnarray}\label{4.10}
\left\{
\begin{array}{ll}
\phi_t-\psi_x=0,\\[2mm]
\psi_t-\big(\frac{P}{V}\phi\big)_x+\big(\frac{R}{V}\xi\big)_x
=Q_1-R_{1,x}, \\[2mm]
\frac{R}{\gamma-1}\xi_t+P \psi_x
+\left(\frac{R\xi}{V+\phi}-\frac{P\phi}{V+\phi}\right)U_x
-\big(\frac{\kappa }{V}\xi_x\big)_x+
\big(\frac{\kappa \Theta_x}{V^2}\phi\big)_x
=Q_2+UR_{1,x}-R_{2,x}.
\end{array}
\right.
\end{eqnarray}
where
\begin{eqnarray*}
Q_1:&=&-\bigg(\frac{R\xi}{V+\phi}-\frac{R}{V}\xi-\frac{P\phi}{V+\phi}
+\frac{P}{V}\phi\bigg)_x\\[2mm]
&=&\bigg(\frac{R\xi}{(V+\phi)^2}-
\frac{P(2V+\phi)\phi}{V(V+\phi)^2}
\bigg)\phi_x+
O(1)(|\phi\xi_x|+|\phi\xi|+|\phi|^2)
\\[2mm]
&=&O(1)(|\phi\phi_x|+|\xi\phi_x|+|\phi\xi_x|+|\phi\xi|+|\phi|^2), \nonumber\\[2mm]
Q_2:&=&-\bigg(\frac{R\xi}{V+\phi}-\frac{P\phi}{V+\phi}\bigg)\psi_x
-\bigg(\frac{\kappa \Theta_x}{vV^2}\phi^2\bigg)_x+\bigg(\frac{\kappa \phi \xi_x}{vV}\bigg)_x\\[2mm]
&=&O(1)(|\phi\psi_x| +|\xi \psi_x|
+|\phi|^2 + |\phi\xi_x|
+|\phi \phi_x|+|\phi_x\xi_x|+|\phi \xi_{xx}|).\nonumber
\end{eqnarray*}
\begin{Lemma}It holds that
\begin{eqnarray}
&&\|(\phi,\psi,\xi)(t)\|^2_1+\int^t_0 \|(\xi_x,\xi_{xx})(\tau)\|^2 \,d\tau
\\[2mm]
&\leq &C(\|(\phi_0,\psi_0,\xi_0)\|^2_1+\delta^{\frac{1}{2}}+|\beta_2|)+
C(\epsilon_0+\delta_0)\int^t_0{\|(\phi,\psi,\xi)(\tau)\|^2_2}d\tau. \nonumber
\end{eqnarray}
\end{Lemma}
{\bf Proof.}\ \ Multiplying $(\ref{4.10})_1$ by $\phi$, $(\ref{4.10})_2$ by $\frac{V}{P}\psi$ and
$(\ref{4.10})_3$ by $ \frac{R}{P^2}\xi$ respectively and adding all the resultant equations, we
have
\begin{eqnarray}\label{4.13}
&&\bigg\{\frac{\phi^2}{2}+\frac{V}{2P}\psi^2+\frac{R^2}{2(\gamma-1)P^2}\xi^2\bigg\}_t
+\bigg\{-\phi\psi+\frac{R}{P}\psi \xi+\frac{R\kappa\Theta_x}{V^2P^2}\phi\xi-
\frac{R\kappa}{VP^2}\xi\xi_x\bigg\}_x\nonumber\\[2mm]
&&+J_3+J_4=\frac{V}{P}\psi( Q_1- R_{1,x})+ \frac{R}{P^2}\xi( Q_2+UR_{1,x}-R_{2,x}),
\end{eqnarray}
where
\begin{eqnarray*}
J_3:&=&-\big(\frac{V}{2P}\big)_t \psi^2+\frac{R\kappa}{VP^2}\xi^2_x+
\frac{\kappa}{V}\big(\frac{R}{P^2}\big)_x\xi\xi_x+
\frac{R^2P_t}{(\gamma-1)P^3}\xi^2,\\[2mm]
J_4:&=&\frac{R}{V}\big(\frac{V}{P}\big)_x\psi \xi
-\frac{V}{P}\big(\frac{P}{V}\big)_x\phi \psi
-\big(\frac{R}{P^2}\xi\big)_x\frac{\kappa \Theta_x}{V^2}\phi
+\frac{R}{P^2}\xi\big(\frac{R\xi}{V+\phi}-\frac{P\phi}{V+\phi}\big)U_x.
\end{eqnarray*}
Noting that the right hand side of (4.14) can be treated similarly as in the proof of
Lemma 4.1,
and integrating (4.14) on ${\bf{R}}\times (0,t)$, we have
\begin{eqnarray}\label{4.14}
&&\|(\phi,\psi,\xi)(t)\|^2+\int^t_0 \|\xi_x(\tau)\|^2\,d\tau
\\[2mm]
&\leq &C(\|(\phi_0,\psi_0,\xi_0)\|^2+\delta^{\frac{1}{2}}+|\beta_2|)
+C(\epsilon_0+\delta_0)\int^t_0{\|(\phi,\psi,\xi)(\tau)\|^2_2}\,d\tau.\nonumber
\end{eqnarray}
Similarly, multiplying $(\ref{4.10})_{1,x}$ by $\phi_x$, $(\ref{4.10})_{2,x}$ by $\frac{V}{P}\psi_x$ and
$(\ref{4.10})_{3,x}$ by $ \frac{R}{P^2}\xi_x$ respectively, adding them all and
integrating the resultant formula on ${\bf{R}}\times(0,t)$, we obtain
\begin{eqnarray}\label{4.16}
&&\|(\phi_x,\psi_x,\xi_x)(t)\|^2+\int^t_0\|\xi_{xx}(\tau)\|^2\,d\tau\\[2mm]
&\leq &C(\|(\phi_0,\psi_0,\xi_0)\|^2_1+\delta^{\frac{1}{2}}+|\beta_2|)+
C(\epsilon_0+\delta_0)\int^t_0{\|(\phi,\psi,\xi)(\tau)\|^2_2}\,d\tau,\nonumber
\end{eqnarray}
where we used the integration by parts
$$
\int \frac{R}{P^2}\xi_x Q_{2,x}\,dx = -\int \big(\frac{R}{P^2}\xi_x\big)_x Q_{2}\,dx.
$$
Combing (4.15) and (4.16), we can have the desired estimate (4.13).
This completes the proof of Lemma 4.2.\\

\begin{Lemma}If $\delta_0$ and $\epsilon_0$ are suitably small, it holds that
\begin{eqnarray}\label{4.19}
&&\|(\phi_{xx},\psi_{xx},\xi_{xx})(t)\|^2+\int^t_0
\|\xi_{xxx}(\tau)\|^2\,d\tau\\[2mm]
&&\leq C(\|(\phi_0,\psi_0,\xi_0)\|^2_2+\delta^{\frac{1}{2}}+|\beta_2|)
+C(\epsilon_0+\delta_0)\int^t_0{\|(\phi,\psi,\xi)(\tau)\|^2_2}\,d\tau.\nonumber
\end{eqnarray}
\end{Lemma}
{\bf Proof.}\ \ Multiplying $(\ref{4.10})_{1,xx}$ by $\phi_{xx}$, $(\ref{4.10})_{2,xx}$ by $\frac{V}{P}\psi_{xx}$ and
$(\ref{4.10})_{3,xx}$ by $ \frac{R}{P^2}\xi_{xx}$ respectively and adding them all,
then integrating the resultant formula on ${\bf{R}}\times [0,t]$,  we can get with the aid of
(4.13) that
\begin{eqnarray}\label{4.21}
&&\|(\phi_{xx},\psi_{xx},\xi_{xx})(t)\|^2
+\int^t_0 \|\xi_{xxx}(\tau)\|^2\,d\tau \nonumber\\[2mm]
&\leq &C(\|(\phi_0,\psi_0,\xi_0)\|_2^2+\delta^\frac{1}{2}+|\beta_2|)+
C(\epsilon_0+\delta_0)\int^t_0{\|(\xi,\phi,\psi)(\tau)\|^2_2}\,d\tau\\[2mm]
&&+C|\int^t_0\int{\frac{V}{P}\psi_{xx}Q_{1,xx}}\,dx\,d\tau|+
C\int^t_0|\int
\big(\frac{R}{P^2}\xi_{xx}\big)_xQ_{2,x}\,dx|\,d\tau,\nonumber
\end{eqnarray}
where we should emphasize that we have to treat the third term of the left hand side of
(4.18) because $Q_{1}$ includes a principal part of the quasi-linear hyperbolic system for
$(\phi,\psi)$.
First, the last term in the right hand side of (4.18) is easily estimated by
the Sobolev's inequality as
\begin{equation}
\int^t_0|\int
\big(\frac{R}{P^2}\xi_{xx}\big)_xQ_{2,x}\,dx|\,d\tau
\le C(\epsilon_0+\delta_0)\int^t_0(\|(\phi,\psi,\xi)(\tau)\|^2_2+\|\xi_{xxx}\|^2)\,d\tau.
\end{equation}
As for the third term, recalling
$$
Q_1=A(\phi,\xi,V)\phi_x+
O(1)(|\phi\xi_x|+|\phi\xi|+|\phi|^2)
$$
where
$$
A(\phi,\xi,V) = \frac{R\xi}{(V+\phi)^2}-
\frac{P(2V+\phi)\phi}{V(V+\phi)^2},
$$
we have by integration by parts
\begin{eqnarray}\label{4.23}
&&\int^t_0\int {\frac{V}{P}\psi_{xx}Q_{1,xx}}\,dxd\tau\nonumber\\[2mm]
&&=-\int^t_0\int \frac{V}{P}\psi_{xxx}A\phi_{xx}\,dxd\tau
-\int^t_0\int \big(\frac{V}{P}A\big)_x\psi_{xx}\phi_{xx}\,dxd\tau
\\[2mm]
&&\quad
+\int^t_0\int \frac{V}{P}\psi_{xx}\big(2A_x\phi_{xx}+A_{xx}\phi_x+
\big(O(1)(|\phi\xi_x|+|\phi\xi|+|\phi|^2)\big)_{xx}\big)\,dxd\tau\nonumber\\[2mm]
&&=: I_1+I_2+I_3.\nonumber
\end{eqnarray}
By the relation $\phi_t=\psi_x$, we further rewrite $I_1$ as
\begin{eqnarray*}
I_1&=& -\int^t_0\int \frac{V}{P}A\phi_{xx}\phi_{xxt}\,dxd\tau\\[2mm]
&=&  -\frac{1}{2}\int^t_0\int \big(\frac{V}{P}A|\phi_{xx}|^2\big)_t\,dxd\tau
+  \frac{1}{2}\int^t_0\int \big(\frac{V}{P}A\big)_t|\phi_{xx}|^2\,dxd\tau  \\[2mm]
&=& I_{11}+I_{12}.
\end{eqnarray*}
We estimate the integrals one by one as follows.
\begin{eqnarray}
|I_{11}| &\le&  \frac{1}{2}|\int^t_0\int \big(\frac{V}{P}A|\phi_{xx}|^2\big)_t\,dxd\tau|\nonumber\\[2mm]
&\le& \frac{1}{2}\bigg|\left[\int \frac{V}{P}A|\phi_{xx}|^2\,dx\right]^t_0\bigg|\\[2mm]
&\le& C(\epsilon_0\|\phi_{xx}(t)\|^2+\|\phi_{0,xx}\|^2);\nonumber
\end{eqnarray}
\begin{eqnarray}
|I_{12}| &\le& \frac{1}{2}\int^t_0\int
\big|(\frac{V}{P})_tA+ \frac{V}{P}A_VV_t+\frac{V}{P}A_\phi\phi_t+\frac{V}{P}A_\xi\xi_t\big|
|\phi_{xx}|^2\,dxd\tau \nonumber \\[-1mm]
&&{\ }\\[-1mm]
&\le& C(\epsilon_0+\delta_0)\int^t_0(\|\phi_{xx}(\tau)\|^2+\|\xi_{xx}(\tau)\|_1^2)\,d\tau\nonumber.
\end{eqnarray}
In this process to have the estimate (4.22), the term
$$
\int^t_0 \int |\xi_{xx}||\phi_{xx}|^2\,dxd\tau
$$
which appears after inserting the equation of $\xi$ into $\xi_t$
should be carefully estimated  by using the Sobolev's inequality as
\begin{eqnarray*}
\int^t_0\int |\xi_{xx}||\phi_{xx}|^2\,dxd\tau &\le&
\int^t_0\|\xi_{xx}\|_1\|\phi_{xx}\|^2\,d\tau \\[2mm]
&\le& \epsilon_0
\int^t_0(\|\xi_{xx}(\tau)\|_1^2+\|\phi_{xx}(\tau)\|^2)\,d\tau.
\end{eqnarray*}
Similarly, noting that $O(1)$ appeared in $I_3$ is regarded as smooth functions of
$(\phi,\psi,\xi,V)$, we can show
\begin{equation}
|I_2|+|I_3| \le C(\epsilon_0+\delta_0)
\int^t_0(\|(\phi,\psi,\xi)(\tau)\|^2_2+\|\xi_{xxx}(\tau)\|^2)\,d\tau.
\end{equation}
Therefore, inserting the estimates (4.19)-(4.23) into (4.18) and
choosing $\epsilon_0$ and $\delta_0$ suitably small, we have the
desired estimate (4.17).
This completes the proof of Lemma 4.3.\\

Combing  up the results of Lemma 4.1 to Lemma 4.3, and choosing
$\epsilon_0$ and $\delta_0$ suitably small,  we have
\begin{eqnarray}\label{4.24}
&&\|(\Phi,\Psi,W)(t)\|^2
+\|(\phi,\psi,\xi)(t)\|^2_2\nonumber
\\[2mm]
&&\quad
+\int^t_0\int \big((|U_{1x}|+|U_{3x}|)(\Psi^2+W^2)
+W_{x}^2\big)\,dxd\tau
+\int^t_0{\|\xi(\tau)\|^2_3}\,d\tau\\[2mm]
&\leq &C(\|(\phi_0,\psi_0,\xi_0)\|^2_2+\delta^{\frac{1}{2}}+|\beta_2|)
+C(\epsilon_0+\delta_0)\int^t_0{\|(\phi,\psi)(\tau)\|^2_2}\,d\tau.\nonumber
\end{eqnarray}
As for the last term of the right hand of (\ref{4.24}), we show
the following lemma which is very essential in this paper, in contrast with
the previous paper [2].
\begin{Lemma}
If $\delta_0$ and $\epsilon_0$ are suitably small, it holds that
\begin{eqnarray}\label{4.25}
\int^t_0{\|(\phi,\psi)(\tau)\|^2_2}\,d\tau
\leq C(\|(\phi_0,\psi_0,\xi_0)\|^2_2+\delta^{\frac{1}{2}}+|\beta_2|).
\end{eqnarray}
\end{Lemma}
{\bf Proof.}\ \ Multiplying the equation $(4.1)_2$ by $-\frac{P}{2}\Phi_x$,
$(4.1)_3$ by $\Psi_x$ respectively and adding the resultant equations, we
have
\begin{eqnarray}\label{4.26}
&&
\bigg\{\frac{R}{(\gamma-1)}W\Psi_x-\frac{P}{2}\Phi_x\Psi\bigg\}_t
+\bigg\{\frac{P}{2}\Phi_t \Psi -\frac{R}{(\gamma-1)}W\Psi_t\bigg\}_x
+\frac{P^2}{2V}\Phi^2_x+\frac{P}{2}\Psi^2_x+J_5\nonumber\\[2mm]
&=&\big(\frac{P}{2} \Phi_x+\frac{R}{(\gamma-1)}W_x\big)( N_1- R_1)
+ \Psi_x\big(\frac{\kappa}{V}W_{xx}+ N_2-R_2\big),
\end{eqnarray}
where
\begin{eqnarray*}
J_5=O(1)\bigg(|\Phi_x W_x|+|W_x|^2
+|U_x|(|\Phi_x \Psi|+|\Psi_x\Psi|+|\Phi_x\Psi_x|+ |\Psi_x|^2+|W_x\Psi|)\bigg).
\end{eqnarray*}
Integrating (\ref{4.26}) and choosing $\delta_0$ suitably small, we have
\begin{eqnarray}\label{4.27}
\int^t_0\int{(\Phi^2_x+\Psi^2_x)}\,dxd\tau
&\leq&C(\|(\Psi,W,\Phi_x ,\Psi_x)(t)\|^2+\|(\Psi_0,W_0,\Phi_{0x} ,\Psi_{0x})\|^2)\nonumber\\[2mm]
&&+\, C\int^t_0\int{((|U_{1,x}|+|U_{3,x}|)\Psi^2+W^2_x+W^2_{xx})}\,dxd\tau\\[2mm]
&&+\, C\int^t_0\int{((|\Phi_x|+|W_x|)|N_1- R_1|+|\Psi_x||N_2-R_2|+|\tilde{R}|\Psi^2)}\,dxd\tau.\nonumber
\end{eqnarray}
Estimating the last term of (\ref{4.27}) as in the proof of Lemma 4.1 and also using (\ref{4.24}),  we can obtain
\begin{eqnarray}\label{4.28}
\int^t_0\int {(\Phi^2_x+\Psi^2_x)}\,dxd\tau
\leq C(\|(\phi_0,\psi_0,\xi_0)\|^2_2+\delta^{\frac{1}{2}}+|\beta_2|)
+C(\epsilon_0+\delta_0)\int^t_0{\|(\phi,\psi)(\tau)\|^2_2}\,d\tau.
\end{eqnarray}
Similarly,
multiplying $(\ref{4.10})_2$ by $-\frac{P}{2}\phi_x$,
$(\ref{4.10})_3$ by $\psi_x$ and multiplying $(\ref{4.10})_{2,x}$ by $-\frac{P}{2}\phi_{xx}$,
$(\ref{4.10})_{3,x}$ by $\psi_{xx}$ respectively,
adding them all and integrating the resultant formula,
we can also get by using (\ref{4.24}) and (\ref{4.28})
\begin{eqnarray}\label{4.32}
&&\int^t_0\int{(\phi^2_x+\psi^2_x+\phi^2_{xx}+\psi^2_{xx})}\,dxd\tau\\[2mm]
&\leq& C(\|(\phi_0,\psi_0,\xi_0)\|^2_2+\delta^{\frac{1}{2}}+|\beta_2|)
+C(\epsilon_0+\delta_0)\int^t_0{\|(\phi,\psi)(\tau)\|^2_2}\,d\tau,\nonumber
\end{eqnarray}
Putting (\ref{4.28}) and (\ref{4.32}) together
and taking $\delta_0$ and $\epsilon_0$  suitably small,
we have the desired estimate (\ref{4.25}) immediately. This completes the proof of Lemma 4.4.\\

Inserting  (\ref{4.25}) into (\ref{4.24}) and recalling the relations
$$
\phi=\Phi_x,\quad \psi=\Psi_x, \quad
 \xi=W_x+\frac{\gamma-1}{R}(U_x\Psi-\frac{1}{2}\Psi^2_x),
$$
we finally reach at the desired {\it a priori} estimates for suitably small
$\delta_0$ and $\epsilon_0$:
\begin{eqnarray*}
&&\|(\Phi,\Psi,W)(t)\|^2_3+\int^t_0\int(|U_{1x}|+|U_{3x}|)(\Psi^2+W^2)(x,\tau)\,
dxd\tau\\[2mm]
&&+\int^t_0{(\|(\Phi_x,\Psi_x,W_x)(\tau)\|^2_2+\|\xi(\tau)\|^2_3)}\,d\tau
\leq C_0(\|(\Phi_0,\Psi_0,W_0)\|^2_3+\delta^{\frac{1}{2}}+|\beta_2|).
\end{eqnarray*}
Thus the proof of Proposition 3.2 is completed.


\begin{thebibliography}{99}

\bibitem{Huang-Li-Matsumura} F.-M. Huang, J. Li  and  A. Matsumura.
Asymptotic stability of combination of viscous contact wave with rarefaction waves
for one-dimensional compressible Navier-Stokes system.
Arch. Rational Mech. Anal. {\bf 197} (2010), 89-116.

\bibitem{Huang-Matsumura} F.-M. Huang  and A. Matsumura,  Stability of a composite wave of two
viscous shock waves  for the full compressible Navier-Stokes equation.
{\it Comm. Math. Phys.} {\bf 289} (2009), 841-861.

\bibitem{Huang-Matsumura-1}  F.-M. Huang, A. Matsumura and  Z.-P. Xin,
Stability of contact discontinuities for the 1-D
compressible Navier-Stokes equations.
Arch. Ration. Mech. Anal. {\bf 179} (2006), 55-77.

\bibitem{Huang-Zhao} F.-M.  Huang and H.-J. Zhao,
On the global stability of contact discontinuity
for compressible Navier-Stokes equations.
Rend. Sem. Mat. Univ. {\bf 109} (2003), 283-305.

\bibitem{Goodman} J. Goodman, Nonlinear asymptotic stability
of viscous shock profiles for conservation laws.
Arch. Rational Mech. Anal. {\bf 95}(1986), 325-344

\bibitem{kawashima} S. Kawashima,
Large-time behavior of solutions to hyperbolic-parabolic
systems of conservation laws and applications,
Proc. Roy. Soc. Edinburgh Sect. A  106  (1987), 169-194.

\bibitem{kawashima-thesis} S. Kawashima, Systems of a Hyperbolic-Parabolic Composite Type, with Applications to the Equations of Magnetohydrodynamics,  doctoral thesis in Kyoto University, 1984.
    http://hdl.handle.net/2433/97887


\bibitem{Kawashima-Matsumura} S. Kawashima and A. Matsumura, Asymptotic stability of
travelling wave solutions of systems for one-dimensional gas motion.
{\it Comm. Math. Phys.} {\bf 101} (1985), 97-127.

\bibitem{T.P.Liu-1} T.-P. Liu, Shock wave for  compresible Navier-Stokes equations are stable.
{\it Comm. Math. Phys.} {\bf 50} (1986), 565-594.

\bibitem{T.P.Liu-2}  T.-P. Liu, Nonlinear stability of shock waves for viscous conservation laws.
Mem. Am. Math. Soc. {\bf 56}(1985), 1-108.


\bibitem{Matsumura-Nishihara-1}   A. Matsumura and K. Nishihara, Asymptotics toward the rarefaction
waves of a one-dimensional model system for compressible viscous gas.
Japan J. Appl. Math. {\bf 3}(1985), 3-13

\bibitem{Matsumura-Nishihara-2}  A. Matsumura and K. Nishihara, Global stability
of the rarefaction wave of a one-dimensional model system for compressible viscous gas.
Comm. Math. Phys. {\bf 144}(1992), 325-335.

\bibitem{Matsumura-Nishihara-4} A. Matsumura and  K. Nishihara, Global stability toward the rerafaction
wave for the solutions of  viscous p-system with  boundary effect,
 Quart. Appl. Math.  {\bf 58} (2000), 69-83.



\bibitem{Nishihara-Yang-Zhao} K. Nishihara, T. Yang and  H.-J. Zhao,
Nonlinear stability of strong rarefaction waves for compressible
Navier-Stokes equations. SIAM J. Math. Anal. {\bf 35} (2004), 1561-1597.

\bibitem{Murakami} T. Murakami,
Asymptotics toward the rarefaction
waves for a one-dimensional model system of
non-viscous and heat-conductive ideal gas, in preparation.



\bibitem{Smoller} J. Smoller, {\it Shock wave and Reaction-Diffusion Equations}, Second Edition,
New York: Springer-Verlag, 1994

\bibitem{Szepessy} A. Szepessy and Z.-P. Xin, Nonlinear stability of viscous shock waves.
Arch. Ration. Mech. Anal. {\bf 122}(1993), 53-103.

\end{thebibliography}
\end{document}